\documentclass[reqno]{amsart}
\usepackage{amsmath}
\usepackage{amsthm}
\usepackage{amssymb}
\usepackage{graphics}

\usepackage[latin1]{inputenc} 

\theoremstyle{definition}
\newtheorem{teorema}{Theorem}
\newtheorem*{constants}{Constants}

\newtheorem*{notation}{Notation}
\newtheorem{propo}{Proposition}[section]
\newtheorem{lema}{Lemma}[section]
\newtheorem*{notacion}{Notation}
\newtheorem*{nota}{Remark}
\newtheorem{notanum}{Remark}[section]
\newtheorem*{defi}{Definition}

\newcommand{\R}{\mathbb R}

\newcommand{\Z}{\mathbb Z}
\newcommand{\N}{\mathbb N}

\newcommand{\s}{\sigma}
\renewcommand{\b}{\beta}

\renewcommand{\d}{\delta}
\renewcommand{\l}{\lambda}
\renewcommand{\t}{\theta}
\renewcommand{\a}{\alpha}

\newcommand{\D}{\mathcal D}

\newcommand{\q}{\hat{q}}

\renewcommand{\t}{\theta}

\newcommand{\suma}{\displaystyle\sum}
\newcommand{\union}{\displaystyle\bigcup}

\newcommand{\norm}[1]{\| #1 \|}
\newcommand{\biggnorm}[1]{\biggl\| #1 \biggr\|}

\newcommand{\tocho}{\q(\xi)\q(\eta-\tau)\q(\tau-\xi)}

\newcommand{\e}{|\eta|}
\newcommand{\es}{|\eta(s)|}
\newcommand{\xim}{|\xi|}
\newcommand{\taum}{|\tau|}
\newcommand{\dif}{|\eta-\xi|}

\newcommand{\esca}{\xi\cdot (\eta-\xi)}
\newcommand{\escal}{\tau\cdot (\eta-\tau)}

\newcommand{\difmt}{\Big|\tau-{\eta\over 2}\Big|}

\newcommand{\Gae}{\Gamma (\eta )}
\newcommand{\Lax}{\Lambda (\xi)}

\newcommand{\Gamase}{\Gamma_{\epsi}^+ (\eta )}
\newcommand{\Gamenose}{\Gamma_{\epsi}^- (\eta )}

\newcommand{\Gainf}{\Gamma_{\infty} (\eta )}
\newcommand{\Gaj}{\Gamma_{j} (\eta )}

\newcommand{\homoa}{\dot{W}^{\a,2}}
\newcommand{\homoame}{\dot{W}^{\a-1+\epsi,2}}
\newcommand{\homos}{\dot{W}^{s,2}}

\newcommand{\homomenos}{\dot{W}^{-{1\over 2},2}}

\newcommand{\f}[2]{\frac{#1}{#2}}
\newcommand{\Rn}{\mathbb{R}^{n}}
\newcommand{\ra}{\rightarrow}
\newcommand{\Rdos}{\mathbb{R}^{2}}

\newcommand{\Gad}{\Gamma_{\d}(\eta)}

\allowdisplaybreaks

\newcommand {\epsi} {\varepsilon}

\newcommand {\La} {{\Lambda}}

\thanks{This paper was published in \emph{Inverse Problems} \textbf{23} (2007) 625-643; doi:10.1088/0266-5611/23/2/010}

\begin{document}

\title{Inverse backscattering for the Schrödinger equation in 2D}
\author{{ Juan Manuel Reyes}}

\maketitle

\begin{center}
{\small Department of Mathematics and Statistics\\ University of Helsinki\\
  FI-00014 Helsinki, Finland\\
E-mail: reyes.juanmanuel@gmail.com }
\end{center}

\begin{abstract}
We study the inverse backscattering problem for the Schrödinger
equation in two dimensions. We prove that, for a non-smooth
potential  in 2D the main singularities up to $1/2$  of the
derivative of the potential are contained in the Born
approximation (Diffraction Tomography  approximation) constructed
from the backscattering data. We measure singularities in the
scale of Hilbertian Sobolev spaces.
\end{abstract}

\section{Introduction}

We consider the inverse scattering problem for the Schrödinger
operator $H=-\Delta + q(x)\,$, with a  real-valued potential
$q(x)$. The scattering solution $u=u(k,\t , x)$ associated with
the energy $k^2 $ and the incident direction $\t$ is defined as
the solution of the problem
\begin{equation}\label{1.1}
\left\{
\begin{array}{ll}
(-\Delta + q -k^2 )u = 0\, , & \\
u = e^{ikx\cdot \,\t } + u_{s}\, , & \\
\end{array}
\right.
\end{equation}
where the function $u_{s}$ satisfies the outgoing Sommerfeld
radiation condition, which means, for a compactly supported
potential $q\,$, that $u$ has asymptotics as $|x|\ra +\infty$
\begin{equation}
\label{asymptotics}u(k , \t , x) = e^{ik x\cdot \,\t} + C
|x|^{{1-n\over 2}} k^{{n-3\over 2}} e^{ik|x|} A(k , \t , \t' ) +
o(|x|^{{1-n\over 2}})\, ,
\end{equation}
where $\t' = {x\over |x|}\,$. The function $A(k , \t , \t' )\,$,
$x\in\R\,$, $\t\, ,\, \t' $ in the unit sphere $S^{n-1}\,$, is
called the scattering amplitude or far-field pattern.

The  outgoing resolvent   operator for  the  Laplacian  is given,
in  terms of the  Fourier  transform,  by
\[
\widehat{R_{k} (f)}(\xi)=  (-\xim^2 + k^2 + i0 )^{-1}
\widehat{f}(\xi)\, .
\]
We obtain the  so-called Lippmann-Schwinger integral equation by
applying the outgoing resolvent to (\ref{1.1})
\begin{equation}
\label{lippschwinger}u_{s}(k , \t , x )=R_{k} (q(\cdot ) e^{ik
(\cdot )\cdot \,\t })(x) + R_{k}(q(\cdot )u_{s}(k , \t , \cdot
))(x)\, .
\end{equation}
The key operator in the above integral equation is
\[
T_{k}(f)(x) = R_{k}(q(\cdot )f(\cdot ))(x)\, .
\]
There are several {\it a priori} estimates for $R_{k}$ that allow
us to prove existence and uniqueness of Lippmann-Schwinger
integral equation. Usually,  Fredholm theory applies and
everything follows from compactness arguments, the Rellich
uniqueness theorem and unique continuation principles, in the case
of real-valued potentials. The solution can be obtained in several
situations (these cases do not require $q$ to be  real) by
perturbation arguments, assuming that the energy is sufficiently
large, $k>k_{0}\geq 0\,$, where $k_{0}$ depends on some {\it a
priori} bound of the potential $q\,$. As an example, we may
consider compactly supported $q\in L^r (\Rn)$ for some $r>{n\over
2}\,$. In this case, the resolvent operator $R_k $ is bounded from
$L^p (\Rn)$ to $L^{p' }(\Rn )$ with the norm decaying to $0$ as
$k\ra \infty $ when ${1\over p} - {1\over p' } = {1\over r}\,$,
see \cite{A}, \cite{KRS} and see also \cite{R1}. This together
with Hölder inequality proves that for big $k$ the operator
$T_{k}$ is a contraction in $L^p $ and then existence and
uniqueness of solution of \eqref{lippschwinger} easily follow.

Once the scattering solution is obtained we may prove that the
far-field pattern can be expressed as
\begin{equation}
\label{farfield}A(k , \t , \t' )=\int_{\Rn} e^{-ik\t' \cdot \, y}
q(y) u(k , \t , y) dy\, ,
\end{equation}
see \cite{ER2} where this  is  used  as  a definition  for  the
non-compactly supported $q$.

The Born series of $q$
is obtained by inserting the Lippmann-Schwinger integral equation
in \eqref{farfield}
\begin{equation}
\label{bornn}A(k , \t , \t') = \q(k(\t' - \t)) +
\suma_{j=2}^{\infty}\widehat{Q_j^{\t } (q)}(k(\t' - \t ))\, ,
\end{equation}
where
\[
\widehat{Q_j^{\t} (q)}(k(\t' - \t )) = \int_{\Rn} e^{-ik\t'
\,\cdot\, y} (q R_k )^{j-1} (q(\cdot )e^{ik\t\cdot
(\cdot)})(y)dy\, .
\]

We deal with the backscattering inverse problem, for which one
assumes the data with the direction of the receiver opposed to the
incident  direction (echoes),  i.e $A(k, \t , -\t )\,$. The
inverse problem is  then formally well determined. The unique
determination of  $q$ by  the   backscattering  data is  an open
problem. Local  and  generic  uniqueness  have  been  proved  by
Eskin and Ralston \cite{ER2},  see  also \cite{S}.

In this  case, we obtain the Born series for backscattering data
\begin{equation}
\label{bornseries}
A(k , \t , -\t ) = \q(\xi) +
\suma_{j=2}^{\infty} \widehat{Q_{j} (q) }(\xi)\, ,
\end{equation}
where $\xi = -2k\t $ and the $j$-adic term in the Born series is
given by
\[
\widehat{Q_j (q)}(\xi) = \int_{\Rn} e^{ik\t\,\cdot\, y} (q R_k
)^{j-1} (q(\cdot ) e^{ik\t\cdot (\cdot)})(y)dy\, .
\]
We  define  the Born  approximation  for  backscattering  data  as
\[
\widehat{q_B }(\xi)=A(k, \t , -\t )
\]
where $\xi=-2k\t\,$.

The   approximated   potential  $q_B$  is  the  target of Diffraction  Tomography.
  In  this  paper,  we  study  how  much information  on  the  actual  potential $q$  can  be  obtained  from $q_B$. We are  able  to  prove  that  the  main  singularities  of  $q$ and $q_B$   are  the  same   up  to $1/2^-$  derivative.

A procedure  to  obtain  this  recovery of  singularities from
Diffraction Tomography  is  to  give  regularity estimates  for
the  $j$-adic  term  in the Born  series (\ref{bornseries}). The
first  of  these  estimates  in 2D  was obtained  by  Ola,
P\"aiv\"arinta  and Serov, see \cite{OPS}, concerning  the
quadratic term $Q_2$  and  was improved  by Ruiz  and Vargas, see
\cite{RV},  who obtained  the  mentioned  $1/2^-$ gaining of
derivative for  the quadratic  term.   Nevertheless  for  $q$ in the  Sobolev
space $ W^{s,2}$ with $s>1/2$,  the  known estimates  for the term
$Q_3$ in  the Born series are not  sufficient  to  assure  that
$q-q_B \in W^{\alpha,2}$ for  $\alpha<s+1/2$.  To  achieve  this
we prove the  main  result  of  this  work.

\begin{teorema}
\label{theoremone}{\it Assume that $q$ is a compactly supported
function in $W^{s\, ,\, 2}(\Rdos)\,$, for $0\leq s<1\,$. Then
$Q_{3}(q)\in W^{\a\, ,\, 2}(\Rdos) + C^{\infty}(\Rdos)\,$, for any
$0\leq\a<s+1\,$.}

\end{teorema}

This  theorem, together  with  the Ruiz  and Vargas   estimate
for the quadratic  term  and  their  estimates  for  the  general
$j$-adic term with  $j>3$,  allows  us  to  claim

\begin{teorema}
\label{theoremtwo}{\it Assume that $q$ is a compactly supported
function in $W^{s\, ,\, 2} \, $, for $0\leq s <1\,$. Then
$q-q_{B}\in W^{\a\, ,\, 2}+C^{\infty}\,$, for any $\a\in\R$ such
that $0\leq \a<s+{1\over 2}\,$.}

\end{teorema}

We expect each term in the Born series (\ref{bornseries}) to win
half a derivative with respect to the previous one, claiming the
conjecture $Q_j (q)\in W^{\a_j , 2}(\Rdos) + C^{\infty}(\Rdos)\,$,
for all $\a_j \in\R $ with $0\leq \a_j <s + {j-1\over 2}\quad
(j\geq 2)\,$, provided that $q$ is a compactly supported function
and $q\in W^{s , 2}(\Rdos)\,$, $0\leq s<1\,$. We address this
question in a future work.

The  results in  this paper and in \cite{RV} could  be  extended
for  non-compactly  supported  potential  assuming  an appropriate
decay  at infinity. To  simplify  the matter  we reduce ourselves
to  the  compactly  supported  case.

Section 2 is the main one in this paper. We prove theorem
\ref{theoremtwo} in Section 3. In Section 4, we include some
lemmata often used in Section 2. In particular, lemma
\ref{lemmasum} is
essential in order to get estimate \eqref{estimatetwo}.\\

\begin{constants}
We use the letter $C$ to denote any constant that can be
explicitly computed in terms of known quantities. The exact value
denoted by $C$ may therefore change from line to line in a given
computation.
\end{constants}

\begin{notation}
We will use the following notation for the Hilbertian Sobolev
space and the homogeneous Hilbertian Sobolev space, respectively:
\begin{align*}
&W^{s,\, 2} = \{ f\in\D'(\Rn) : (1+|\cdot|^2 )^{{s\over 2}}
\widehat{f}(\cdot)\in L^2 \}\, ,\\
&\dot{W}^{s,\, 2} = \{ f\in\D'(\Rn) : D^s f =
\mathcal{F}^{-1}\left(\, |\cdot|^s \widehat{f}(\cdot)\, \right
)\in L^2 \}\, .
\end{align*}

The expression $|\xi-\tau|\sim 2^{-k}\e$ means
that $ 2^{-k-1}\e\leq |\xi-\tau|\leq 2^{-k+1}\e\, $. In this work
$\chi$ denotes the characteristic function of the set
$\{\eta\in\Rdos : \e>10\}\,$. The letter $M$ denotes the
Hardy-Littlewood
 maximal operator. We denote the one-dimensional
 Hausdorff measure in $\Rdos\,$ by $\s\,$. Let $\eta,\, \xi\in\Rdos\setminus\{0\}\,$. We write
\begin{equation}
\label{notation1}\Gae :=\left\{ x \in\Rdos : \left| x -{\eta\over
2}\right|={\e\over 2}\right\}\, ,
\end{equation}
referring to the circumference centered at ${\eta\over 2}$ and
radius ${\e\over 2}\,$ and
\begin{equation}
\label{notation2}\Lambda(\xi):=\left\{ x \in\Rdos : \xi\cdot ( x
-\xi )=0 \right\}
\end{equation}
denotes the line orthogonal to $\xi$ that contains the point $\xi\,$.\\
\end{notation}

\section{Proof of theorem \ref{theoremone}}

The cubic term in the Born series for backscattering data is given
by
\[
\widehat{Q_{3}(q)}(\xi):=\int_{\Rdos}e^{ik\t\cdot y}(qR_{+}(k^2
))^2 (q(\cdot)e^{ik\t\cdot (\cdot )})(y)dy\, ,
\]
for any $\xi\in\Rdos\,$, where $\xi=-2k\t\,$, that is,
$k={\xim\over 2}$ and $\t=-{\xi\over \xim}\,$. From lemma 3.1 in
\cite{R2}, this term can be characterized by the following:

\begin{propo}
\label{prop3.2}{\it For any $\eta\in\Rdos\setminus\{0\}\,$,
\begin{align}
\label{pvterm}\widehat{Q_{3}(q)}(\eta)&=
p.v.\int_{\Rdos}\int_{\Rdos}{\q(\xi)\q(\eta-\tau)\q(\tau-\xi)\over
[\esca]\, [\tau\cdot (\eta-\tau)]}\, d\xi d\tau\\
\notag&+ 2{i\pi\over \e}\,
p.v.\int_{\Rdos}\int_{\Gae}{\q(\xi)\q(\eta-\tau)\q(\tau-\xi)\over
\tau\cdot (\eta-\tau)}\, d\s(\xi)d\tau\\
\notag&- {\pi^2 \over \e^2
}\int_{\Gae}\int_{\Gae}\q(\xi)\q(\eta-\tau)\q(\tau-\xi)\,
d\s(\tau)d\s(\xi)\, .
\end{align}}
\end{propo}

\begin{notacion}
For any $\eta\in\Rdos\setminus\{0\}\,$, we use the following
notation:
\begin{align}
\notag&\widehat{Q'(q)}(\eta):={1\over \e^2
}\int_{\Gae}\int_{\Gae}\q(\xi)\q(\eta-\tau)\q(\tau-\xi)\,
d\s(\tau)d\s(\xi)\, ,\quad \text{and}\\
\label{exp12}&\widehat{Q''(q)}(\eta):={1\over \e }\, p.v.
\int_{\Rdos}\int_{\Gae}{\q(\xi)\q(\eta-\tau)\q(\tau-\xi)\over
\tau\cdot (\eta-\tau)}\, d\s(\xi)d\tau\, .
\end{align}
\end{notacion}


In fact, we are going to prove that
\[
Q'(q)\, , Q''(q)\in W^{\a,2}(\Rdos)+C^{\infty}(\Rdos)\, ,
\]
for any $\a\,$ such that $0\leq\a<s+1\,$.

\subsection{Estimate of $Q'(q)$}$ $\\

Let us split the set $\Gae\times\Gae$ into the two regions:
\begin{align*}
&I(\eta):=\left\{(\xi , \tau )\in\Gae\times\Gae : |\xi-\tau|\geq
{\e\over 100}\right\}\quad\text{and}\\
&II(\eta):=\left\{(\xi , \tau )\in\Gae\times\Gae : |\xi-\tau|<
{\e\over 100}\right\}\, .
\end{align*}
In this way, we can write $Q'(q)=Q'_{I}(q)+Q'_{II}(q)\,$, where
\[
\widehat{Q'_{I}(q)}(\eta):={1\over\e^2
}\int\int_{I(\eta)}\q(\xi)\q(\eta-\tau)\q(\tau-\xi)\,
d\s(\tau)d\s(\xi)\, ,
\]
for any $\eta\in\Rdos\setminus\{0\}\,$, and an analogous
definition for $\widehat{Q'_{II}(q)}(\eta)\,$. We will prove that
\begin{align}
\label{estimateone}&\left\|\mathcal{F}^{-1}\left(\chi\widehat{Q'_{II}(q)}\right)\right\|_{\homoa}\leq
C\|q\|_{L^2 }^2 \|q\|_{\dot{W}^{\a-1+\epsi,2}}\, , \quad \text{and}\\
\label{estimatetwo}&\left\|\mathcal{F}^{-1}\left(\chi\widehat{Q_{I}'(q)}\right)\right\|_{\homoa}\leq
C(\epsi )\, \left(\|q\|_{L^2 }\|q\|_{\homomenos} + \|q\|_{L^2 }^2
\right)\,\|q\|_{\dot{W}^{\a-1+\epsi\,,\, 2}}\, ,
\end{align}
if $\epsi >0\,$, $0<\a + \epsi <2\,$, where $C(\epsi)$ is a
positive constant depending on $\epsi\,$.

{\it Proof of estimate \eqref{estimateone}.} We know that
$II(\eta)\subset II_{<}(\eta)\cup II_{>}(\eta)\,$, where
\begin{align*}
&II_{<}(\eta):=\left\{(\xi,\tau)\in\Gae\times\Gae :
|\xi-\tau|<{\e\over 100}\, , \, \xim, |\tau|\leq
\left({\sqrt{2}\over 2}+{1\over 100}\right)\e \right\}\, ,\quad
\text{and}\\
&II_{>}(\eta):=\left\{(\xi,\tau)\in\Gae\times\Gae :
|\xi-\tau|<{\e\over 100}\, , \, \dif , |\eta-\tau|\leq
\left({\sqrt{2}\over 2}+{1\over 100}\right)\e \right\}\, .
\end{align*}

By taking the change of variables $\xi'=\eta-\xi\,$,
$\tau'=\eta-\tau\,$, by Fubini's theorem and the symmetry property
$(\xi,\tau)\in II_{>}(\eta)\, \Leftrightarrow \, (\tau,\xi)\in
II_{>}(\eta)\,$, we have
\begin{equation}
\label{simetriaII}\left|\widehat{Q'_{II}(q)}(\eta)\right| \leq 2\,
\widehat{Q'_{II_{>}}(q)}(\eta)\, ,
\end{equation}
where $\widehat{Q'_{II_{>}}(q)}(\eta):={1\over\e^2
}\int\int_{II_{>}(\eta)}|\q(\xi)\q(\eta-\tau)\q(\tau-\xi)|\,
d\s(\tau)d\s(\xi)\,$. Applying the Cauchy-Schwartz inequality to
the last expression and by the properties of the region $II_{>}\,$
we may write
\begin{align}
\notag&\left\|\mathcal{F}^{-1}\left(\chi\widehat{Q'_{II_{>}}(q)}\right)\right\|_{\homoa}^2
= \int_{\{\eta : \e\geq
10\}}\e^{2\a}\left|\widehat{Q'_{II_{>}}(q)}(\eta)\right|^2
d\eta\\
\label{exp13}&\leq\int_{X}\e^{2\a-4} \int_{\{\xi\in\,\Gae \, :\,
\xim\geq \, C_{2}\e \}}|\q(\xi)|^2 \int_{\Gae}|\q(\eta-\tau)|^2
d\s(\tau)\, d\s(\xi)\, F(\eta)\, d\eta\, ,
\end{align}
where $X:=\{\eta\in\Rdos :\, \e>10 \}\,$ ,
$F(\eta)=\int\int_{\Gae\times\Gae}|\q(\tau'-\xi')|^2 d\s(\tau')
d\s(\xi')\,$ and $C_{2}:=\left(1-\left({1\over\sqrt{2}}+{1\over
100}\right)^2 \right)^{{1\over 2}}\,$. We know that
\begin{align*}
&\{(\eta , \xi )\in\Rdos\times\Rdos : \e>10\, ,\, \xi\in\Gae\, ,\,
\xim\geq C_{2}\e \}\\
&\subset\{(\eta,\xi)\in\Rdos\times\Rdos : \xim>1,\, \e\leq
C_{2}^{-1}\xim\, ,\, \eta\in\Lax\}\, ,
\end{align*}
and by lemma \ref{lema2.2}, we may change the order of integration
and estimate expression \eqref{exp13} by
\begin{align}
\notag&\leq \int_{\{\xi\in\Rdos : \xim>1\}}{|\q(\xi)|^2 \over
|\xi|}\int_{\Omega(\xi)
}\e^{2\a-3}\int_{\Gamma(\eta)}|\q(\eta-\tau)|^2 d\s(\tau)\, F(\eta)\, d\s(\eta)\, d\xi\\
\label{exp14}&\leq C \, \norm{q}_{L^2 }^2 \int_{\{\xi\in\Rdos :
\xim >1\}}{|\q(\xi)|^2 \over
|\xi|}\int_{\Omega(\xi)}\e^{2\a-2}\int_{\Gamma(\eta)}|\q(\eta-\tau)|^2
d\s(\tau)\, d\s(\eta)\, d\xi\, ,
\end{align}
where $\,\Omega (\xi):=\{\eta\in\Lax : \e\leq \, C_{2}^{-1}\xim
\}\,$, and we have used the inequality $F(\eta)\leq C
\e\norm{q}_{L^2 }^2 \, $. Let us see the last inequality. If we
widen the curve $\Gae$ until
$\Gamma_{1}(\eta):=\left\{\tau\in\Rdos : \left|\, \difmt-{\e\over
2} \,\right|<1\right\}\,$, by part $1)$ of lemma \ref{remark} we
have:
\[
F(\eta)\leq C\int_{\Gae}\int_{\Gamma_{1}(\eta)}M\q(\tau-\xi')^2
d\tau d\s(\xi')\leq C\norm{M\q}_{L^2 }^2 \,\s(\Gae)\leq
C\e\norm{q}_{L^2 }^2 \, ,
\]
where the last inequality follows from the boundedness of
Hardy-Littlewood maximal o\-pe\-ra\-tor
 in
$L^2 (\Rn)\,$ and Plancherel identity, since the measure of $\Gae$
is $\pi\e\,$. In the same way,
\[
\int_{\Gae}|\q(\eta-\tau)|^2 d\s(\tau)\leq
C\int_{\Gamma_{1}(\eta)} \left| M\q(\eta-\tau')\right|^2
d\tau'\leq C \int_{\Rdos}\left|M\q(x)\right|^2 dx\leq
C\norm{q}_{L^2 }^2 \, .
\]
Obviously, if $\eta\in\Omega (\xi)\,$ then $\xim\sim\e\,$.
Expression \eqref{exp14} can be bounded by
\begin{align*}
&C \, \norm{q}_{L^2 }^4 \int_{\{\xi\in\Rdos : \xim
>1\}}{|\q(\xi)|^2 \over
|\xi|}\xim^{2\a-2}\s\left(\Omega(\xi)\right)\, d\xi\leq
C\norm{q}_{L^2 }^4 \norm{q}_{\dot{W}^{\a-1+\epsi\, , \, 2}}^2 \, .
\end{align*}
So we have proved estimate \eqref{estimateone}.

\begin{flushright} $\square $ \end{flushright}

{\it Proof of estimate \eqref{estimatetwo}.} Taking the change of
variable $\tau=\eta-\tau'\,$, we have
\begin{align*}
\widehat{Q'_{I}(q)}(\eta)&={1\over\e^2
}\int\int_{I(\eta)}\q(\xi)\q(\eta-\tau')\q(\tau'-\xi)\,
d\s(\tau')d\s(\xi)\\
&= {1\over\e^2 }\int\int_{\{(\xi,\tau)\in \, \Gae\times\Gae :
|\xi-(\eta-\tau)|\geq {\e\over
100}\}}\q(\xi)\q(\tau)\q(\eta-\tau-\xi)\, d\s(\tau)d\s(\xi)\, .
\end{align*}
Notice that if $\e\geq 4\,$ we can write
\[
\{(\xi,\tau)\in \, \Gae\times\Gae : |\xi-(\eta-\tau)|\geq {\e\over
100}\}=\union_{k=2}^{\left[\log_{2}\e\right]} I_{k}(\eta)\cup
I_{0}(\eta)\cup I_{\infty}(\eta)\, ,
\]
where for any $k\in\mathbb{Z}$ such that $2\leq k\leq
\left[\log_{2}\e\right]\,$ we denote:
\begin{align}
\label{7b}&I_{k}(\eta):=\left\{(\xi,\tau)\in\Gae\times\Gae :
|\xi-\tau|\sim 2^{-k}\e ,\, |\xi-(\eta-\tau)|\geq {\e\over
100}\right\}\, ,\\
\notag&I_{0}(\eta):=\left\{(\xi,\tau)\in\Gae\times\Gae :
|\xi-\tau|\geq {\e\over 2},\, |\xi-(\eta-\tau)|\geq {\e\over 100}
\right\}\, ,\\
\notag&I_{\infty}(\eta):=\left\{(\xi,\tau)\in\Gae\times\Gae :
|\xi-\tau|\leq 1,\, |\xi-(\eta-\tau)|\geq {\e\over 100} \right\}\,
.
\end{align}

Note that cases $k=0$ and $k=\infty\,$ are needed since the union
from $k=2$ to $k=\left[\log_{2}\e\right]\,$ only covers the set of
$\xi\, ,\,\tau\,$ such that $1\leq |\xi-\tau|\leq {\e\over 2}\,$.

For each $k\in\{2,3,...\}\,$ we define
\begin{align*}
&\widehat{Q_{I_{k}}'(q)}(\eta):= \chi_{\{\eta : \e\geq \, 2^k
\}}(\eta)\, {1\over \e^2
}\int\int_{I_{k}(\eta)}|\q(\xi)\q(\tau)\q(\eta-\tau-\xi)|\,
d\s(\tau)d\s(\xi) \,  ,
\end{align*}
and the same expression for $I_{0}(\eta)\, ,\,I_{\infty}(\eta)\,
$, but with $\chi_{\{\eta : \e\geq \, 10\}}(\eta) \,$. For any
$\e\geq 10\,$,
\[
\left|\widehat{Q_{I}'(q)}(\eta)\right|\leq
\suma_{k=2}^{+\infty}\left|\widehat{Q_{I_{k}}'(q)}(\eta)\right| +
\left|\widehat{Q_{I_{0}}'(q)}(\eta)\right| +
\left|\widehat{Q_{I_{\infty}}'(q)}(\eta)\right|\, ,
\]
and then to prove \eqref{estimatetwo} we use
\[
\left\|\mathcal{F}^{-1}\left(\chi\widehat{Q_{I}'(q)}\right)\right\|_{\homoa}\leq
\suma_{k=2}^{+\infty}\|Q_{I_{k}}'(q)\|_{\homoa}+\|Q_{I_{0}}'(q)\|_{\homoa}+\|Q_{I_{\infty}}'(q)\|_{\homoa}\,
.
\]
Let $\epsi\, ,\, \a$ be real numbers with $\epsi >0\,$. For each
$k\in\{2,3,...\}\,$ we claim
\begin{equation}
\label{estimatetwoa}\|Q_{I_{k}}'(q)\|_{\homoa}\leq C\cdot
2^{-\epsi k}\, \|q\|_{L^2
}\|q\|_{\homomenos}\|q\|_{\dot{W}^{\a-1+\epsi,2}}\, .
\end{equation}
Assume $0<\a+\epsi<2\,$. Then we claim
\begin{align}
\label{estimatetwob}&\norm{Q_{I_{\infty}}'(q)}_{\homoa}\leq
C\norm{q}_{L^2
}\left(\norm{q}_{\homomenos}\norm{q}_{\homoame}+\norm{q}_{L^2
}\norm{q}_{\dot{W}^{\a-{3\over 2}+\epsi ,\, 2}}\right)\,
\text{ and}\\
\label{estimatetwoc}&\|Q_{I_{0}}'(q)\|_{\homoa}\leq C \,
\|q\|_{L^2 }\|q\|_{\homomenos}\|q\|_{\dot{W}^{\a-1+\epsi,2}}\, .
\end{align}

In the following, we use the notation in lemma \ref{lemmasum},
which is the key of the proof of the above claims.

\newpage

{\it Proof of claim \eqref{estimatetwoa}.} We take
$I_{k}(\eta)=I_{k}^1 (\eta)\cup I_{k}^2 (\eta)\,$, where $I_{k}^1
(\eta):= \left\{ (\xi,\tau)\in I_{k}(\eta) : |\tau-\eta|\geq
{2^{-k}\e\over 100} \right\}$ and \[ I_{k}^2 (\eta):= \left\{
(\xi,\tau)\in I_{k}(\eta) : |\tau-\eta|\leq {2^{-k}\e\over 100}
\right\}\, .
\]
For each $j\in\{1,2\}\,$, let us define $\widehat{Q_{I_{k}^j
}'(q)}(\eta)$ in the obvious way multiplying by $\chi_{\{\eta :
\e\geq 2^k \}}(\eta )\,$. By Cauchy-Schwartz inequality, and for
$\e\geq 2^k \,$, $j\in\{1,2\}\,$:
\begin{align}
\label{form12}\widehat{Q_{I_{k}^j }'(q)}(\eta) \leq {1\over \e^2
}&\left(\int\int_{I_{k}^j (\eta)}|\q(\xi)\q(\tau)|^2
d\s(\tau)d\s(\xi) \int\int_{I_{k}^j (\eta)}
|\q(\eta-\tau'-\xi')|^2 d\s(\tau')d\s(\xi') \right)^{{1\over 2}}.
\end{align}
Let us begin with $Q_{I_{k}^1 }'(q)\,$. By lemma \ref{remark}, we
have
\begin{align}
\notag\int\int_{I_{k}^1 (\eta)}|\q(\xi)\q(\tau)|^2
d\s(\tau)d\s(\xi) &\leq C \|q\|_{L^2 }^2 \int_{\Psi_{k}(\eta)}
|\q(\tau)|^2 d\s(\tau)\, ,
\end{align}
where $\Psi_{k} (\eta) :=\left\{\tau\in\Gae : |\eta-\tau|\geq
{2^{-k}\e\over 100}\right\}\,$. Since $\e^{-2\epsi}\leq 2^{-2\epsi
k}\,$ and by lemma \ref{lema2.2},
\[
\|Q_{I_{k}^1 }'(q)\|_{\homoa}^2
\]
\begin{align}
\notag&\leq C\cdot 2^{-2\epsi k} \|q\|_{L^2 }^2
\int_{\Rdos}\e^{2\a-4+2\epsi} \int_{\Psi_{k}(\eta)}|\q(\tau)|^2
d\s(\tau) \int\int_{I_{k}(\eta)} |\q(\eta-\tau'-\xi')|^2
d\s(\tau')d\s(\xi')d\eta \\
\notag&= C\cdot 2^{-2\epsi k} \|q\|_{L^2 }^2
\int_{\Rdos}{|\q(\tau)|^2 \over |\tau|}\int_{\Omega_{k} (\tau)
}\e^{2\a-3+2\epsi }\int\int_{I_{k} (\eta)}
|\q(\eta-\tau'-\xi')|^2 d\s(\tau')d\s(\xi') d\s(\eta)d\tau\\
\notag&= C\cdot 2^{-2\epsi k} \|q\|_{L^2 }^2
\int_{\Rdos}{|\q(\tau)|^2 \over |\tau|}F_k (\tau)\, d\tau\leq
C\cdot 2^{-2\epsi k}\,\|q\|_{L^2 }^2 \,\|q\|_{\homomenos}^2
\,\|q\|_{\dot{W}^{\a-1+\epsi , 2}}^2 \, ,
\end{align}
where the last inequality follows from part $i)$ of lemma
\ref{lemmasum} with $C_1 ={1\over 100}$ and $F_k (\tau)\, ,
\Omega_{k}(\tau) \,$ are defined in \eqref{lemma4.4a},
\eqref{omegak}.

We can bound the term $\widehat{Q_{I_{k}^2 }'(q)}(\eta)\,$ in a
similar way. Firstly, we estimate the factor
\[
\int\int_{I_{k}^2 (\eta)}|\q(\xi)\q(\tau)|^2 d\s(\tau) d\s(\xi)
\]
by $C\norm{q}_{L^2 }^2 \int_{\widetilde{\Psi}_{k}(\eta)}
|\q(\xi)|^2 d\s(\xi)\,$, where
$\widetilde{\Psi}_{k}(\eta):=\{\xi\in\Gae : |\eta-\xi|\geq
{49\over 100}\, 2^{-k} \e \}\,$, by using lema \ref{remark} as
above. In order to estimate the expression $\norm{Q_{I_{k}^2
}'(q)}_{\homoa}^2 $ as before, we proceed similarly so that the
variable $\xi$ now acts just as the variable $\tau$ before,
obtaining that
\[
\norm{Q_{I{_k}^2 }'(q)}_{\homoa}^2 \leq C\, 2^{-2\epsi
k}\,\norm{q}_{L^2 }^2 \int_{\Rdos}{|\q(\xi)|^2 \over \xim} \, F_k
(\xi) d\xi\leq C\cdot 2^{-2\epsi k}\,\|q\|_{L^2 }^2
\,\|q\|_{\homomenos}^2 \,\|q\|_{\dot{W}^{\a-1+\epsi , 2}}^2 \,\, ,
\]
where the last inequality follows from part $i)$ of lemma
\ref{lemmasum} with $C_{1}={49\over 100}\,$.

\begin{flushright} $\square $ \end{flushright}

{\it Proof of claim \eqref{estimatetwob}.} We split the set
$I_{\infty}(\eta)$ into $I_{\infty}^1 (\eta)$ (where $\taum\geq
1$) and $I_{\infty}^2 (\eta)$ (where $\taum\leq 1$). Let us denote
$\widehat{Q_{I_{\infty}^j }'(q)}(\eta)$ in the obvious way, for
$j=1,2\,$. In order to estimate the Sobolev norm of
$Q_{I_{\infty}^j }'(q)\, ,\, j=1,2\,$, we apply Cauchy-Schwartz
inequality as in \eqref{form12}, by lemma \ref{remark}, for
$I_{\infty}^1 (\eta)\,$, we may do $|\q(\eta-\xi'-\tau')|\leq C
M\q(\eta-2\tau')\,$, since $|\xi'-\tau'|\leq 1\,$, and for
$I_{\infty}^2 (\eta)\,$, $|\q(\eta-\xi'-\tau')|\leq C
M\q(\eta-\tau')\,$, since $|\xi'|\leq 2\,$, we bound the integral
involving $\q(\xi)$ by lemma \ref{remark}, leading for
$I_{\infty}^1 (\eta)\,$ to $\int\int_{I_{\infty}^1 (\eta)}
M\q(\eta-2\tau' )^2 d\s(\xi')d\s(\tau')\leq
C\int_{\Gae}M\q(\eta-2\tau')^2 d\s(\tau')$ (the same with
$M\q(\eta-\tau')\,$, for $I_{\infty}^2 (\eta)$), change the order
of integration in $\tau\, ,\, \eta$ by lemma \ref{lema2.2}, and
finally, by parts $ii)$ and $iii)$ of lemma \ref{lemmasum}
(provided that $0<\a+\epsi <2$) we get
\begin{align*}
&\norm{Q_{I_{\infty}^1 }'(q)}_{\homoa}^2 \leq C\norm{q}_{L^2 }^2
\left( \norm{q}_{\homomenos}^2 \norm{q}_{\homoame}^2 +
\norm{q}_{L^2 }^2 \norm{q}_{\dot{W}^{\a-{3\over 2}+\epsi \, ,\,
2}}^2 \right)\quad\text{and}\\
&\norm{Q_{I_{\infty}^2 }'(q)}_{\homoa}^2 \leq C \norm{q}_{L^2 }^2
\norm{q}_{\homomenos}^2 \norm{q}_{\homoame}^2 \, ,
\end{align*}
respectively.

\begin{flushright} $\square $ \end{flushright}

{\it Proof of claim \eqref{estimatetwoc}.} We also split the set
$I_0 (\eta)$ into $I_{0}^1 (\eta)$ (where $|\eta-\tau|\geq
{\e\over 4}$) and $I_0^2 (\eta)$ (where $|\eta-\tau|<{\e\over
4}$). Note that on the region $I_0^2 (\eta)\,$, $|\eta-\xi|\geq
{\e\over 4}$ holds. In both cases we apply Cauchy-Schwartz
inequality in the same way as in the previous cases, bound
$|\q(\xi)|^2 $ (for $I_0^1 (\eta)$) or $|\q(\tau)|^2 $ (for $I_0^2
(\eta)$) by the maximal operator by lemma \ref{remark}, change the
order of integration in the variables $\tau\, ,\, \eta \,$, for
$I_0^1 (\eta)$ (in the variables $\xi\, ,\, \eta\,$, for $I_0^2
(\eta)$) by lemma \ref{lema2.2} and finally, by part $i)$ of lemma
\ref{lemmasum}, with $k=1$ and $C_1 ={1\over 2}\,$, we get
\[
\norm{Q_{I_0^j }'(q)}_{\homoa}^2 \leq C\norm{q}_{L^2 }^2
\norm{q}_{\homomenos}^2 \norm{q}_{\homoame}^2 \, , \quad j=1,2\, .
\]

\begin{flushright} $\square $ \end{flushright}

Hence by estimates \eqref{estimatetwoa}, \eqref{estimatetwob} and
\eqref{estimatetwoc} we can write
\[
\left\|\mathcal{F}^{-1}\left(\chi\widehat{Q_{I}'(q)}\right)\right\|_{\homoa}\leq
{C\cdot 2^{-\epsi}\over 2^{\epsi} -1}\,\norm{q}_{L^2
}\left(\norm{q}_{\homomenos}\norm{q}_{\homoame}+\norm{q}_{L^2
}\norm{q}_{\dot{W}^{\a-{3\over 2}+\epsi ,\, 2}}\right)\, ,
\]
and we have proved \eqref{estimatetwo}.

\begin{flushright} $\square $ \end{flushright}

To obtain the non-homogeneous Sobolev norm we proceed as follows.
By lemma \ref{remark}, $q\in W^{-(1-\epsi),\, 2}(\Rdos)\,$ holds
for $0<\epsi<1\,$, and replacing $\a$ by $0\,$ in
\eqref{estimateone} we get that
\[
\mathcal{F}^{-1}\left(\chi\widehat{Q'_{II}(q)}\right)\in
L^{2}(\Rdos)\, .
\]
Note that estimates \eqref{estimatetwoa}, \eqref{estimatetwob} and
\eqref{estimatetwoc} remain true if $\a=0\,$ (assuming that
$0<\epsi<2$ to guarantee the estimate \eqref{estimatetwob}). Hence
$\mathcal{F}^{-1}\left(\chi\widehat{Q'_{I}(q)}\right)\in
L^{2}(\Rdos)\,$. It follows that
\[
Q'(q)=\mathcal{F}^{-1}\left(\left(1-\chi\right)\widehat{Q'(q)}\right)
+ \mathcal{F}^{-1}\left(\chi\widehat{Q'(q)}\right)\, ,
\]
where the first term is a function belonging to the class
$C^{\infty}(\Rdos)\,$, and the second one is in
$W^{\a,2}(\Rdos)\,$ if $0\leq\a\leq s+1-\epsi\;$ (with $0<\epsi<2$
arbitrary, provided that $s<1\,$), that is to say, if
$0\leq\a<s+1\,$. So, we have finished
with the term $Q'(q)\,$.\\

\subsection{Estimate of $Q''(q)$}$ $\\

The singularities of the integral \eqref{exp12} are those points
$\tau$ in the plane such that $\tau\cdot (\eta-\tau)=0$, that is,
the set $\Gae$. So, we decompose the plane in an annulus
containing $\Gae$ and its complement. Next, we decompose the first
annulus in diadic coronas and try to treat the corresponding
integral terms. Let
\[
N:=\max\{ [\log_{2}\e ]-2,1 \}=\left\{
\begin{array}{ll}
[\log_{2}\e]-2, &\text{if }\e\geq 16,\\
1, &\text{if }\e<16.
\end{array}
\right.
\]
Let $j_{0}$ be the lowest integer such that $j_{0}\geq
-1-\log_{2}(\d_{0})\,$, with $\d_0 $ from lemma \ref{lemadelta}
(see the appendix). We define the sets
\begin{align}
\notag&\Gamma_{j_0^{\, -} }(\eta):=\{ \tau\in\Rdos : \Big|\,
\Big|\tau-{\eta\over 2}\Big|-{|\eta|\over 2}
\, \Big|>2^{-j_0 -1}\e \},\\
\notag&\Gaj:=\{ \tau\in\Rdos : 2^{-j-2}\e < \Big|\,
\Big|\tau-{\eta\over
2}\Big|-{|\eta|\over 2} \, \Big| \leq 2^{-j-1}\e \},\\
\label{gainf}&\Gainf:=\{ \tau\in\Rdos : \Big|\,
\Big|\tau-{\eta\over 2}\Big|-{|\eta|\over 2} \, \Big|\leq 2\},
\end{align}
with $j_0 \leq j\leq N$. If $j\geq j_0 $ it is true that $j\leq N
\, \Leftrightarrow\, \e\geq 2^{j+2}\,$ (for $\e\geq 16 $). So, we
also define for $j_0 \leq j < \infty\,$:
\begin{equation}
\label{qsegundaj}\widehat{Q''_{j}(q)}(\eta):=
\widetilde{\chi}(\e){1\over\e
}\int_{\Gaj}\int_{\Gae}{\q(\xi)\q(\eta-\tau)\q(\tau-\xi)\over
\tau\cdot
(\eta-\tau)}\, d\s(\xi)d\tau\, ,\\
\end{equation}
where $\widetilde{\chi}=\chi_{[2^{j+2},+\infty)}\,$, and the
obvious notations for $\widehat{Q''_{\infty}(q)}(\eta)\,$,
$\widehat{Q''_{j_0^{\, -}}(q)}(\eta)\,$ without the characteristic
function. Since $\Rdos = \bigcup_{j=j_0}^N \Gaj
\cup\Gainf\cup\Gamma_{j_0^{\, -}}(\eta)\,$, for any
$\eta\in\Rdos\setminus\{0\}\,$, it follows that
\[
\widehat{Q''(q)}(\eta)=\widehat{Q''_{j_0^{\, -}}(q)}(\eta) +
\suma_{j=j_0 }^{N}\widehat{Q''_{j}(q)}(\eta)+
\widehat{Q''_{\infty}(q)}(\eta) = \widehat{Q''_{j_0^{\,
-}}(q)}(\eta) + \suma_{j=j_0
}^{\infty}\widehat{Q''_{j}(q)}(\eta)+\widehat{Q''_{\infty}(q)}(\eta)\,
.
\]

It is easy to see that
\begin{equation}
\label{exp100}\biggnorm{\mathcal{F}^{-1}\left(\chi
\,\widehat{Q_{j_0^{\, -}}''(q)}\right)}_{\homoa}\leq
C\norm{q}_{L^2 }^2 \|q\|_{\homoame}
 \, .
\end{equation}\\

{\it Bound of the corona terms.}

By Minkowski's inequality, we have
\begin{align}
\notag&\Big|\Big|\mathcal{F}^{-1}\Big(\chi \suma_{j=j_0
}^{\infty}\widehat{Q''_{j}(q)}\Big)\Big|\Big|_{\homoa}
=\Big|\Big|\suma_{j=j_0 }^{\infty}\mathcal{F}^{-1}\Big( \chi
\widehat{Q''_{j}(q)}\Big)\Big|\Big|_{\homoa}\leq \suma_{j=j_0
}^{\infty}\Big|\Big|\mathcal{F}^{-1}\Big(\chi
\widehat{Q''_{j}(q)}\Big)\Big|\Big|_{\homoa}.
\end{align}
If $j\geq j_0 $ and $\tau\in\Gaj$, $|\tau\cdot (\eta-\tau)|\geq
2^{-j-3}\e^2 \,$, from where we deduce that
\begin{align*}
|\widehat{Q''_{j}(q)}(\eta)| &\leq
2^{j+3}\chi_{(2^{j+1},+\infty)}\Big(\f{\e}{2}\Big){1\over \e^3
}\int_{\Gaj}
\int_{\Gae}\left|\q(\xi)\q(\eta-\tau)\q(\tau-\xi)\right|\,
d\s(\xi)d\tau\, .
\end{align*}

It follows that
\[
\Gaj \subset \{ \tau\in\Rdos :
\Big||\tau-\f{\eta}{2}|-|\f{\eta}{2}|\Big|<2^{-j-1}\e \},
\]
hence, applying the key lemma \ref{lemadelta} with $\d=2^{-j-1}$,
we know that there exist $\d_{0}>0$, $\beta
>1$ and $C>0$ so that for any $j\in\N $ satisfying
$2^{-j-1}<\d_{0}$ (that is, $j\geq j_0 $), we have that
\[
\Big|\Big|\mathcal{F}^{-1}\Big(\chi
\widehat{Q''_{j}(q)}\Big)\Big|\Big|_{\homoa}\leq C
2^{j+3}(2^{-j-1})^{\beta }\|q\|_{\homoame} \left( \norm{q}_{L^2
}\norm{q}_{ \homomenos }+ \norm{q}_{L^2 }^2 \right)\, ,
\]
where $\epsi > 0$ satisfies $0<\a + \epsi <2\,$. We can write
\begin{align}
\notag\suma_{j=j_{0}}^{\infty}\Big|\Big|\mathcal{F}^{-1}\Big(\chi
\widehat{Q''_{j}(q)}\Big)\Big|\Big|_{\homoa} &\leq C
\suma_{j=j_{0}}^{\infty} 2^{-j(\beta-1)} 2^{3-\beta
}\|q\|_{\homoame} \left( \norm{q}_{L^2 }\norm{q}_{ \homomenos
}+ \norm{q}_{L^2 }^2 \right) \\
\label{suma}&= C\|q\|_{\homoame} \left( \norm{q}_{L^2 }\|q\|_{
\homomenos }+ \norm{q}_{L^2 }^2 \right) \, .
\end{align}
The series $\sum_{j=j_{0}}^{\infty}2^{-j(\beta-1)}$ converges
because $\beta>1$.\\

{\it Bound of the singular part close to $\Gae\,$.}

We are going to prove the estimate:
\begin{equation}
\label{exp102}\|\mathcal{F}^{-1}(\chi\,\widehat{Q_{\infty}''(q)})\|_{\homoa}\leq
C \left[\,\norm{q}_{L^2 } \norm{q}_{\homomenos} + \norm{q}_{L^2
}^2 \,\right]\,\norm{q}_{\homoame}\, ,
\end{equation}
for a constant $C>0$ depending on the support of $q\,$, and
provided that $0<\a+\epsi<2\,$. Up to now, we have avoided the
singular region $\Gae\,$. The domain $\Gainf$ contains it. In
order to calculate the principal value of the integral on
$\Gainf$, we integrate on two rings whose radial distance to the
singular circumference is $\epsi>0$ and pass to limit when
$\epsi\ra 0^+ $. We write:
\begin{equation}
\label{exp1}\widehat{Q_{\infty}''(q)}(\eta) =
{1\over\e}\displaystyle\lim_{\epsi\ra 0^+
}\left(\int_{\Gamma_{\epsi}^+ (\eta)}+\int_{\Gamma_{\epsi}^-
(\eta)}\right)\int_{\Gae}\f{\tocho}{\escal}\, d\s(\xi)d\tau\, ,
\end{equation}
where
\begin{align*}
\Gamase &:= \{ \tau\in\Rdos : \epsi < \difmt-\f{\e}{2}<2 \}\text{ and}\\
\Gamenose &:= \{ \tau\in\Rdos : \epsi < \f{\e}{2}-\difmt<2\}.
\end{align*}

Let us take the change of variables $\tau'=\phi(\tau)$,
$\tau\in\Gamenose$, that sends $\tau$ to its symmetrical point
$\tau'\in\Gamase$ with respect to $\Gae$ on the radial direction
with centre at $\f{\eta}{2}\,$. We have
\[
\tau'=\eta-\tau+\e\f{\tau - {\eta\over 2}}{\difmt}.
\]
A straightforward calculation leads up to the following
identities:
\begin{align}
\label{tocho603}&\Big|\phi(\tau)-\f{\eta}{2}\,
\Big|-\f{\e}{2}=-\Big(\difmt-\f{\e}{2}\Big)\, ,\\
\label{tocho604}&|D\phi(\tau)|=1+2\,
\f{\f{\e}{2}-\difmt}{\difmt}\, ,\\
\label{tocho607}&|\phi(\tau)-\tau|=2\Big(\f{\e}{2}-\difmt\Big)\,
,\\
\notag&\phi(\tau)\cdot (\eta-\phi(\tau)) = \Big(\f{\e}{2}+\Big|
\phi(\tau)-\f{\eta}{2}\Big|\Big)\cdot\Big(\difmt-\f{\e}{2}\Big)\, ,\\
\label{tocho602}&\escal =
\Big(\f{\e}{2}+\difmt\Big)\cdot\Big(\f{\e}{2}-\difmt\Big)\, .
\end{align}
Taking the change $\tau'=\phi(\tau)$ in the first integral in
\eqref{exp1}, we get
\begin{align*}
&\int_{\Gamase}\int_{\Gae}\f{\q(\xi)\q(\eta-\tau')\q(\tau'-\xi)}{\tau'\cdot(\eta-\tau')}\,
d\s(\xi) d\tau' \\
&=
\int_{\Gamenose}\int_{\Gae}\f{\q(\xi)\q(\eta-\phi(\tau))\q(\phi(\tau)-\xi)}{\phi(\tau)\cdot(\eta-\phi(\tau))}\,
|D\phi(\tau)|\, d\s(\xi) d\tau\, .
\end{align*}
Then we have   $\quad\widehat{Q_{\infty}''(q)}(\eta)$
\begin{align*}
&=\displaystyle\lim_{\epsi\ra 0^+
}\e^{-1}\int_{\Gamenose}\int_{\Gae}\Big[\,
\f{\q(\xi)\q(\eta-\phi(\tau))\q(\phi(\tau)-\xi)}{\phi(\tau)\cdot
(\eta-\phi(\tau))}\, |D\phi(\tau)|\, + \f{\tocho}{\escal}\,\Big]\,d\s(\xi)d\tau\\
&=\displaystyle\lim_{\epsi\ra 0^+ }\e^{-1}\Big[\,
\int_{\Gamenose}\int_{\Gae}\f{\q(\xi)\q(\eta-\phi(\tau))\q(\phi(\tau)-\xi)-\tocho
}{|\f{\eta}{2}|^2 -\Big|\phi(\tau)-\f{\eta}{2}\Big|^2 }\,
|D\phi(\tau)|\, d\s(\xi)d\tau\\
&+\int_{\Gamenose}\int_{\Gae}\f{\tocho}{|\f{\eta}{2}|^2
 - \Big|\phi(\tau)-\f{\eta}{2}\Big|^2 }\, |D\phi(\tau)| \,
d\s(\xi)d\tau + \int_{\Gamenose}\int_{\Gae}\f{\tocho}{\escal}\,
d\s(\xi)d\tau\, \Big]\\
&=\displaystyle\lim_{\epsi\ra 0^+ }\e^{-1}\left[\,
\int_{\Gamenose}\int_{\Gae}\f{\q(\xi)\left[\q(\eta-\phi(\tau))
\q(\phi(\tau)-\xi)-\q(\eta-\tau)\q(\tau-\xi)\right]}{|\f{\eta}{2}|^2
-\Big|\phi(\tau)-\f{\eta}{2}\Big|^2 }\, |D\phi(\tau)|\,
d\s(\xi)d\tau\right.\\
&-2 \int_{\Gamenose}\int_{\Gae}\f{\tocho}{\Big(|\f{\eta}{2}|
+\Big|\phi(\tau)-\f{\eta}{2}\Big|\Big)\cdot \difmt }\,
d\s(\xi)d\tau\\
&\left.+2 \int_{\Gamenose}\int_{\Gae}\f{ \tocho}{
\Big(|\f{\eta}{2}|+|\phi(\tau)-\f{\eta}{2}|\Big)
\Big(\f{\e}{2}+\difmt\Big)}\; d\s(\xi)d\tau\right]\\
&=:\displaystyle\lim_{\epsi\ra 0^+ }(I_{1}^\epsi + I_{2}^\epsi
+I_{3}^\epsi )\, ,
\end{align*}
where we have to keep identities \eqref{tocho603},
\eqref{tocho604} and \eqref{tocho602} in mind and also
\begin{align*}
&{1\over |\f{\eta}{2}|^2 -\Big|\phi(\tau)-\f{\eta}{2}\Big|^2 }\, +
\f{1}{\escal}=\f{2}{
\Big(|\f{\eta}{2}|+|\phi(\tau)-\f{\eta}{2}|\Big)
\Big(\f{\e}{2}+\difmt\Big)}\, .
\end{align*}

For $\e>10$, the terms $I_{2}^\epsi,\ I_{3}^\epsi$ may be upper
bounded by a term like
\begin{equation}
\label{tocho605}\widehat{J(q)}(\eta):=\chi (\eta)
\int_{\Gamenose}\int_{\Gae}\f{|\tocho|}{\e^3 }\, d\s(\xi)d\tau.
\end{equation}
If one replaces the characteristic function $\chi$ of the set
$\{\eta\in\Rdos : \e >10\}$ by the characteristic function of the
complement of a bigger ball our proof for theorem \ref{theoremone}
remains valid. In \eqref{tocho605} if we replace $\chi$ by the
characteristic function of the set $\{\eta\in\Rdos : \e >r \}\,$,
with $r\in\R$ such that \label{elerre}$r >{2\over \d_0}\,$ (for
$\d_0$ from lemma \ref{lemadelta}) then $\widehat{J(q)}(\eta)\leq
\widehat{Q_{{2\over r }}''(q)}(\eta)\,$ holds, according to the
notation from lemma \ref{lemadelta} (since $2={2\over\e}\,\e
<{2\over r }\, \e\,$ and $\Gainf\subset \Gamma_{{2\over r
}}(\eta)\,$), that is, we may apply lemma \ref{lemadelta} with $\d
= {2\over r} \, (<\d_0)\,$ and get that there exists a constant
$C>0$ such that
\begin{equation}
\label{tocho606}\|J(q)\|_{\homoa}\leq C\left(\,\norm{q}_{L^2
}\|q\|_{\homomenos}\|q\|_{\homoame} + \norm{q}_{L^2 }^2
\norm{q}_{\homoame}\,\right) .
\end{equation}

On the one hand,
\begin{align*}
&I_{1}^\epsi  = {1\over\e
}\int_{\Gamenose}\int_{\Gae}\f{\q(\xi)\left[\,\q(\eta-\phi(\tau))\q(\phi(\tau)-\xi)
-\q(\eta-\tau)\q(\tau-\xi)\,\right]}{|\f{\eta}{2}|^2
-\Big|\phi(\tau)-\f{\eta}{2}\Big|^2 }\, d\s(\xi)d\tau\\
&+ {2\over\e}
\int_{\Gamenose}\int_{\Gae}\f{\q(\xi)\left[\,\q(\eta-\phi(\tau))\q(\phi(\tau)-\xi)
-\q(\eta-\tau)\q(\tau-\xi)\,\right]}{|\f{\eta}{2}|^2
-\Big|\phi(\tau)-\f{\eta}{2}\Big|^2 }\,
\f{\f{\e}{2}-\difmt}{\difmt} \, d\s(\xi)d\tau\, .
\end{align*}
If $\tau\in\Gamenose$ and $\e>10$,
\[
0<\f{\f{\e}{2}-\difmt}{\difmt}<1.
\]
That is,
\begin{align*}
|I_{1}^\epsi| &\leq {3\over\e} \int_{\Gamenose}\int_{\Gae}\f{|
\q(\xi)\q(\phi(\tau)-\xi)\, [\q(\eta-\phi(\tau))-\q(\eta-\tau)]\,
|}{\left|\,|\f{\eta}{2}|^2
-\Big|\phi(\tau)-\f{\eta}{2}\Big|^2 \right|}\, d\s(\xi)d\tau\\
&+ {3\over\e} \int_{\Gamenose}\int_{\Gae}\f{|
\q(\xi)\q(\eta-\tau)\, [\q(\phi(\tau)-\xi)-\q(\tau-\xi)]\,
|}{\left|\,|\f{\eta}{2}|^2 -\Big|\phi(\tau)-\f{\eta}{2}\Big|^2
\right|}\, d\s(\xi)d\tau\\
&=: \widehat{J_{1}(q)}(\eta) + \widehat{J_{2}(q)}(\eta)\, .
\end{align*}
The term $\widehat{J_{1}(q)}(\eta)$ may be bounded by Calderón
estimate (see Section 2 in \cite{H}):
\[
|f(x)-f(y)|\leq C\, (M(\nabla f)(x)+M(\nabla f)(y))\, |x-y|\quad
\text{a.e.,}
\]
provided that $f\in\dot{W}^{1,p}(\Rn):=\{g\in
\mathcal{D}'(\Rn)\,\, / \,\, \nabla g\in L^p (\Rn)\}\,$, for some
$p>1\,$. So, by \eqref{tocho607} we attain that
\begin{align*}
&|\widehat{J_{1}(q)}(\eta)| \leq C\,
{1\over\e}\int_{\Gamenose}\int_{\Gae}\f{[M(\nabla
\q)(\eta-\phi(\tau))+M(\nabla \q)(\eta-\tau)]\,
|\q(\xi)\q(\phi(\tau)-\xi)|}{\Big|\phi(\tau)-\f{\eta}{2}\Big|+|\f{\eta}{2}|
}\,
d\s(\xi)d\tau\\
&=: C\, (\widehat{J_{1}^1 (q)}(\eta)+\widehat{J_{1}^2
(q)}(\eta))\, .
\end{align*}
Let $\hat{f}:=M(\nabla \q)\,$. It holds
\begin{align}
\notag\widehat{J_{1}^1
(q)}(\eta)&={1\over\e}\int_{\Gamase}\int_{\Gae}\f{\hat{f}(\eta-\tau')\,
|\q(\xi)\q(\tau'-\xi)|}{\Big|\tau'-\f{\eta}{2}\Big|+|\f{\eta}{2}|
}\,
|D\phi^{-1}(\tau')|\, d\s(\xi)d\tau'\\
\label{exp612}&\leq
{1\over\e}\left(\int_{\Gamase}\int_{\Gae}\f{\hat{f}(\eta-\tau')\,
|\q(\xi)\q(\tau'-\xi)|}{\Big|\tau'-\f{\eta}{2}\Big|+|\f{\eta}{2}|
}\,
d\s(\xi)d\tau'\right.\\
\label{exp613}&\left. + \, 2
\int_{\Gamase}\int_{\Gae}\f{\hat{f}(\eta-\tau')\,
|\q(\xi)\q(\tau'-\xi)|}{\Big|\tau'-\f{\eta}{2}\Big|+|\f{\eta}{2}|
}\, \f{|\tau'-\f{\eta}{2}|-\f{\e}{2}}{|\tau'-\f{\eta}{2}|}\,
d\s(\xi)d\tau' \right) .
\end{align}
The second integral \eqref{exp613} is bounded by
\[
\widehat{K(q)}(\eta):=\int_{\Gamase}\int_{\Gae}\f{\hat{f}(\eta-\tau')|\q(\xi)\q(\tau'-\xi)|}{\e^3
}\, d\s(\xi)d\tau'\, .
\]
Applying remark to lemma \ref{lemadelta} with $\d={2\over r}\,$
(where $r$ is defined in page \pageref{elerre}), and by lemma
\ref{radial}, we have
\begin{align}
\notag\|\mathcal{F}^{-1}\left(\chi\,
\widehat{K(q)}\right)\|_{\homoa}&\leq
C\|q\|_{\homoame}\,\left[\,\norm{q}_{L^2 }\norm{f}_{L^2 } +
\norm{q}_{\homomenos}\|f\|_{L^2 } +
\norm{f}_{L^2 }\norm{q}_{L^2 }\,\right]\\
\label{tocho608}&\leq C \left[\,\norm{q}_{L^2 }^2
\norm{q}_{\homoame} + \norm{q}_{\homomenos}\norm{q}_{L^2
}\norm{q}_{\homoame}\,\right]\, .
\end{align}
The integral \eqref{exp612} is bounded by a positive constant
multiplied by
\[
\widehat{K'(q)}(\eta):=\int_{\Gainf}\int_{\Gae}\f{\hat{f}(\eta-\tau')|\q(\xi)\q(\tau'-\xi)|}{\e^2
}\, d\s(\xi)d\tau'\, .
\]
By lemma \ref{lemadeltados} and lemma \ref{radial}, the bound
\eqref{tocho608} works for $\|\mathcal{F}^{-1}\left(\chi\,
K'(q)\right)\|_{\homoa}\,$ too. So we have
\begin{equation}
\label{exp610}\left\|\mathcal{F}^{-1}\left(\chi\, \widehat{J_{1}^1
(q)}\right)\right\|_{\homoa}\leq C \left[\,\norm{q}_{L^2 }^2
\norm{q}_{\homoame} + \norm{q}_{\homomenos}\norm{q}_{L^2
}\norm{q}_{\homoame}\,\right]\, .
\end{equation}
It holds
\begin{align*}
\widehat{J_{1}^2 (q)}(\eta)&\leq C \,
\int_{\Gamenose}\int_{\Gae}\f{M(\nabla \q)(\eta-\tau)\, |\q(\xi)|
M\q(\tau-\xi)}{\e^2 }\, d\s(\xi)d\tau\, ,
\end{align*}
since if $\tau\in\Gamenose$, then $|\phi(\tau)-\tau|<4\,$, and by
lemma \ref{remark}, $|\q(\phi(\tau)-\xi)|\leq C M\q(\tau-\xi)\,$,
for a certain constant $C>0\,$. By lemma \ref{lemadeltados} and
lemma \ref{radial}, we have
\begin{align*}
\left\|\mathcal{F}^{-1}\left(\chi\,
\widehat{J^{2}_{1}(q)}\right)\right\|_{\homoa}&\leq C
\left(\,\norm{q}_{L^2 }^2 +\norm{q}_{\homomenos}\norm{q}_{L^2
}\,\right)\norm{q}_{\homoame}\, .
\end{align*}
This last estimate and \eqref{exp610} lead up to
\[
\left\|\mathcal{F}^{-1}\left(\chi\, \widehat{ J_{1}
(q)}\right)\right\|_{\homoa}\leq C \left[\,\norm{q}_{L^2 }^2
\norm{q}_{\homoame} + \norm{q}_{\homomenos}\norm{q}_{L^2
}\norm{q}_{\homoame}\,\right]\, .
\]
In a similar way, by Calderón estimate, lemma \ref{remark} using
that $|\q(\eta-\tau)|\leq C M\q(\eta-\phi(\tau))\,$, remark to
lemma \ref{lemadelta}, lemma \ref{lemadeltados} and lemma
\ref{radial} we may write
\begin{align*}
\|\mathcal{F}^{-1}\left(\chi\,
\widehat{J_{2}(q)}\right)\|_{\homoa}&\leq C \left(\,\norm{q}_{L^2
}^2 +\norm{q}_{\homomenos}\norm{q}_{L^2 }\,\right)
\norm{q}_{\homoame}\, ,
\end{align*}
and conclude that
\begin{align}
\notag
\|\mathcal{F}^{-1}(\chi\,\widehat{Q_{\infty}''(q)})\|_{\homoa}
&\leq C\left[\,\|\mathcal{F}^{-1}\left(\chi\,
\widehat{J_{1}(q)}\right)\|_{\homoa}+
\|\mathcal{F}^{-1}\left(\chi\,
\widehat{J_{2}(q)}\right)\|_{\homoa}+\|J(q)\|_{\homoa}\, \right]\\
\notag&\leq C \left[\,\norm{q}_{L^2 }^2 \norm{q}_{\homoame} +
\norm{q}_{\homomenos}\norm{q}_{L^2 }\norm{q}_{\homoame}\,\right]\,
,
\end{align}
provided that $0<\a+\epsi<2\,$.\\

It is true that $Q''(q) =
\mathcal{F}^{-1}\left(\left[1-\chi\right]\,
\widehat{Q''(q)}\right) + \mathcal{F}^{-1}\left(\chi\,
\widehat{Q''(q)}\right)\,$, where the first term is a function
belonging to the class $C^{\infty}(\Rdos)\,$, and the second one
belongs to $\dot{W}^{\a\, ,\, 2}(\Rdos)$ taking $\epsi = s+1-\a$
in \eqref{exp100}, \eqref{suma} and \eqref{exp102}. By lemma
\ref{remark}, $q\in \dot{W}^{-1+\epsi\, ,\, 2}(\Rdos)$ (for
$0<\epsi <1$) and replacing $\a$ by zero in \eqref{exp100},
\eqref{suma}, \eqref{exp102} we get that
$\mathcal{F}^{-1}\left(\chi\, \widehat{Q''(q)}\right)\in L^2
(\Rdos)\,$.

The remaining principal value term in \eqref{pvterm} can be
treated in a similar way. We do not include the proof to avoid
tedious repetitions. To control those terms close to singularities
we proceed similarly as we did for $Q_{\infty}''(q)\,$,
compensating signs and using estimates for second differences. We
have finished the proof of theorem \ref{theoremone}.

\begin{flushright} $\square $ \end{flushright}

\section{Proof of theorem \ref{theoremtwo}}

In order to avoid the control of the remainder term in the Born
series, the following proposition gives, modulo a $C^{\infty }$
function, the convergence of the Born series in $W^{\a , 2}\,$,
for $\a < s + {1\over 2}\,$.

\begin{propo}
\label{jestimate}{\it Let $q\in W^{s , 2}$ be a real-valued
compactly supported function for $0\leq s <1\,$. Assume that $C_0
>1\,$. Then, for any $\a\in\R$ such that $\a<s+{1\over 2}\,$:
\begin{equation}
\label{triang}\norm{\widetilde{Q}_j (q)}_{\homoa}\leq C(s , \a
)\,\norm{q}_{L^2 }\norm{q}_{\homos}^{-1} \, A(s, q, j)\, ,
\end{equation}
where $\widetilde{Q}_j (q) = \mathcal{F}^{-1}\left(
\chi^{\ast}\widehat{Q_j (q)} \right)\,$, $\chi^{\ast} (\xi)=0$ if
$\xim\leq C_0\,$, $\chi^{\ast }(\xi)=1$ if $\xim > C_0\,$, $j\geq
4$ and
\[
A(s , q , j ):=\left\{
\begin{array}{ll}
C_0^{{5\over 2}} \left[ 2^{{7\over 8 }} C_0^{-{3\over 4}}
\norm{q}_{\homos}
\right]^j \, , & \text{if } 0\leq s\leq {1\over 2}\, ,\\
C_0^{{25\over 8}} \left[ 2 C_0^{-{7\over 8}} \norm{q}_{\homos}
\right]^j \, , & \text{if } {1\over 2}< s<1\, .

\end{array}
\right.
\]}

\end{propo}

{\it Proof of proposition \ref{jestimate}.} We follow the lines of
proposition 4.3 in \cite{RV}. We lose some regularity in return
for the gain of decay as a negative power of $C_0 \,$. We write
$R_{\t }(k^2 )(f)(x):=e^{-ik\t\cdot\, x} R_+ (k^2 )(e^{ik\t\cdot
(\cdot)} f(\cdot))(x)\,$. By Cauchy-Schwartz inequality:
\begin{equation}
\label{35}\norm{(q R_{\t}(k^2 ))^{j-1}(q)}_{L^1 }\leq
\norm{q}_{L^2 } \norm{R_{\t }(k^2 )(q R_{\t } (k^2
))^{j-2}(q)}_{L^2 }\, ,
\end{equation}
and applying successively the estimate for the resolvent given by
lemma 3.4 in \cite{R1} and the following inequality for Sobolev
spaces due to Zolesio (see \cite{G} and also Section 3.5 in
\cite{T}):
\[
\norm{u\, v}_{W^{s_3 , p}(\Rn)}\leq \norm{u}_{W^{s_1 , p_1 }(\Rn
)}\norm{v}_{W^{s_2 , p_2 }(\Rn)}\, ,
\]
where $s_1 , s_2 , s_3 \geq 0\,$, $s_3 \leq s_1 \,$, $s_3 \leq s_2
\,$, $s_1 +s_2 -s_3 \geq n \left( {1\over p_1 }+{1\over p_2
}-{1\over p } \right)\geq 0$ and $p_j > p\,$, $j=1,2\,$, one can
prove that
\begin{equation}
\label{36}\norm{R_{\t}(k^2 )(q R_{\t }(k^2 ))^{j-2}(q)}_{L^2 }\leq
C k^{-1-\a_j }\norm{q}_{\homos }^{j-1}\, ,
\end{equation}
where
\[
\a_j := \left\{
\begin{array}{ll}
{3\over 4} (j-2) + {s\over 4}(j-1)\, , & \text{if $s\leq {1\over
2}\,$,}\\
(j-3){3+ s\over 4} + 1\, , & \text{if ${1\over 2}\leq s <1\,$.}
\end{array}
\right.
\]
For all $h_j \in\R$ such that $h_j < \a_j $ it holds
\begin{align}
\label{expr3}\norm{\widetilde{Q}_j (q)}_{\dot{W}^{h_j , \, 2}}^2
&\leq C 2^{2 h_j }\int_{{C_0 \over 2}}^{+\infty}k^{2 h_j
+1}\int_{S^1 }\norm{(q R_{\t}(k^2
))^{j-1}(q)}_{L^1 (\Rdos)}^2 d\s(\t) dk\\
\label{expr4}&\leq C 2^{2 h_j }\int_{{C_0 \over 2}}^{+\infty}k^{2
(h_j - \a_j ) -1 } dk\,\, \norm{q}_{L^2 }^2 \norm{q}_{\dot{W}^{s ,
2}}^{2j-2}\, ,
\end{align}
where we pass from \eqref{expr3} to \eqref{expr4} by formulae
\eqref{35} and \eqref{36}, and the last integral in $k$ is
convergent because of $h_j <\a_j \,$. So, we get
\begin{equation}
\label{j} \norm{\widetilde{Q}_j (q)}_{\dot{W}^{h_j , \, 2}}\leq C
{2^{h_j }\over \sqrt{\a_j -h_j }}\, \left({C_0 \over
2}\right)^{h_j -\a_j } \norm{q}_{L^2 }\norm{q}_{\dot{W}^{s ,
2}}^{j-1}\, .
\end{equation}
Let $\epsi=\epsi (s , \a ):= \left(s +{1\over 2}\right)-\a >0\,$.
We have
\begin{align*}
\norm{\widetilde{Q}_j (q)}_{\homoa} &= \left[ \int_{\{\xi : \xim
>C_0 \}}\xim^{2\a }\left| \widehat{Q_j (q)}(\xi)\right|^2 d\xi\right]^{{1\over 2}}\leq
C_0 ^{s +{1\over 2}-\a_j }\norm{\widetilde{Q}_j
(q)}_{\dot{W}^{\a_j -\epsi , 2}} \\
&\leq C {2^{\a_j -\epsi }\over \sqrt{\epsi } } \left({C_0 \over
2}\right)^{-\epsi} C_0 ^{s +{1\over 2}-\a_j }\norm{q}_{L^2 }
\norm{q}_{\dot{W}^{s ,2}}^{j-1} = C(s , \a ) \, 2^{\a_j } \,
C_0^{\a -\a_j } \norm{q}_{L^2 }\norm{q}_{\homos}^{j-1} \, ,
\end{align*}
where last inequality follows from formula \eqref{j} in the case
$h_j = \a_j -\epsi\,$. Since
\[
\a -\a_j < s+{1\over 2} - \a_j \leq \left\{
\begin{array}{ll}
-{3\over 4} j + {5\over 2}\, , & \text{if $s\leq {1\over 2}\,$,}\\
-{7\over 8} j + {25\over 8}\, , & \text{if ${1\over 2} < s
<1\,$,}\\
\end{array}
\right.
\]
and $2^{\a_j } \leq 2^{{7\over 8}j }\,$, if $s\leq {1\over 2}\,$
and $2^{\a_j }\leq 2^j\,$, if ${1\over 2}<s<1\,$, we obtain
\eqref{triang}.

\begin{flushright} $\square $ \end{flushright}

With the notation from proposition \ref{jestimate}, the Born
series \eqref{bornseries} allows us to write
\[
q_B -q = \mathcal{F}^{-1}\left(\chi^{\ast} \,\widehat{q_B
-q}\right) +
\mathcal{F}^{-1}\left(\left(1-\chi^{\ast}\right)\widehat{q_B
-q}\right) = \suma_{j=2}^{\infty} \widetilde{Q}_j (q) +
\mathcal{F}^{-1}\left(\left(1-\chi^{\ast}\right)\, \widehat{q_B -
q}\right)\, ,
\]
where
$\mathcal{F}^{-1}\left(\left(1-\chi^{\ast}\right)\widehat{q_B
-q}\right)$ is $C^{\infty}\,$. If we choose $C_0 $ large enough,
for example, taking
\[
C_0 := \max\{ 10 , 2^{{7\over 6}}\norm{q}_{\homos}^{{4\over 3}} ,
2^{{8\over 7}}\norm{q}_{\homos}^{{8\over 7}} \} + 1\, ,
\]
it is true that $\sum_{j=4}^{+\infty}A(s , q , j)<+\infty \,$.
From theorem \ref{theoremone} and \cite{RV} we can write
\begin{align}
\label{A}&\norm{\widetilde{Q}_2 (q)}_{\homoa}\leq C
\norm{q}_{\homomenos}
\norm{q}_{\homos} \quad \text{and}\\
\label{B}&\norm{\widetilde{Q}_3 (q)}_{\homoa}\leq C \left(
\norm{q}_{L^2 }^2  + \norm{q}_{L^2 }\norm{q}_{\homomenos } \right)
\norm{q}_{\homos}\, ,
\end{align}
for all $\a < s+{1\over 2}\,$, and we have proved that
\begin{equation}
\label{C}\suma_{j=4}^{+\infty}\norm{\widetilde{Q}_j
(q)}_{\homoa}\leq C(s , \a ) C_0^{-{1\over 2}}{ \norm{q}_{L^2 }
\norm{q}_{\homos}^{3} \over 1-2^{{7\over 8}} C_0^{-{3\over 4}}
\norm{q}_{\homos} }\, ,
\end{equation}
if $0\leq s \leq {1\over 2}$ (for ${1\over 2}<s<1\,$, an analogous
expression holds). We know that \eqref{A},\eqref{B},\eqref{C}
remain true if $\a=0\,$, obtaining the non-homogeneous Sobolev
norm. We have finished the proof of theorem \ref{theoremtwo}.

\section*{Acknowledgments}

The work was supported by Grant MTM2005-07652-C02-01 of MEC
(Spain). I would like to express my thanks to my advisor Professor
Alberto Ruiz for his valuable teachings.

\numberwithin{equation}{section}

\section{Appendix}

Let it be the following submanifold of $\mathbb{R}^{2n}\,$ $V:=\{
(\eta,\xi) \in \Rn \times \Rn : \xi\cdot (\xi-\eta)=0 \}\,$. Then
$V$ can be considered from the point of view of the following
spherical sections:
\[
V=\{ (\eta,\xi) \in \Rn \times \Rn : \xi\in\Gae \}\, ,
\]
or the plane sections: $V=\{ (\eta,\xi) \in \Rn \times \Rn :
\eta\in\Lax \}\,$, where $\Gae$ and $\Lax$ are defined in
\eqref{notation1} and \eqref{notation2}. In this context, the
following lemma from \cite{RV} allows us to change the order of
integration in $\xi$ and $\eta\,$.
\begin{lema}
\label{lema2.2}{\it Let $V\equiv \{ (\eta,\xi) \in \Rn \times \Rn
: \xi\cdot (\xi-\eta)=0 \}$. Let $d\s_{\eta}(\xi)$ be the measure
on $\Gae$ induced by the $n$-dimensional Lebesgue measure $d\xi$
and let $d\s_{\xi}(\eta)$ be the measure on $\Lambda(\xi)$ induced
by the $n$-dimensional Lebesgue measure $d\eta$. Then
\[
d\s_{\eta}(\xi)d\eta={\e\over \xim}\, d\s_{\xi}(\eta)d\xi.
\]}
\end{lema}

The following lemma in \cite{RV} is used several times in this
work.
\begin{lema}
\label{remark}{\it Assume that the support of $q$ is contained in
the unit ball. Then we have

$1)$ If $\,\xi, \ \xi' \in \Rn$ satisfy $|\xi-\xi'|\leq 3$, then $
|\q(\xi)|\leq CM\q(\xi')\,$.

$2)$ $\| \q \|_{L^{\infty}}\leq C\| \q\|_{L^{2}}$.

$3)$ For $0<\beta<{n\over 2}$ and $s\in\R\,$, $ \| q
\|_{W^{s-\beta,\, 2}}\leq C \| q \|_{W^{s,\, 2}}\, $, where $C$ depends
 on the size of the support of $q$.}
\end{lema}

We want to indicate a
\begin{defi}
{\it Let $1\leq p<+\infty$. We define \emph{the weights class
$A_{p}$} as the set of the non-negative locally integrable
functions $w$ that satisfy the so-called \emph{condition
$A_{p}\,$}, that is, that there exists a constant $C>0$
independent of $x$ and $r$ so that
\[
\f{1}{|B|^p }\int_{B}w(x)dx\, \left(\int_{B}w(x)^{-{1\over p-1}}\,
dx\right)^{p-1}\leq C
\]
for all ball $B$ centred at $x\in\Rn$ and radius $r>0\,$.}
\end{defi}

Indeed, next lemma \ref{radial}, which is useful to bound the term
$Q_{\infty}''(q)\,$, follows from estimates of the
Hardy-Littlewood maximal operator and checking up on the function
$|x|^{2s}$ belongs to the weighted class $A_{2}$ in two dimensions
if $-1<s<1\,$.

\begin{lema}
\label{radial}{\it Let $q$ be a compactly supported function in
$\dot{W}^{s\, ,\, 2}(\Rdos)\,$. Hence there exists a positive
constant $C$ depending on the support of $q$ such that for any
$s\in\R$ with $|s|<1\,$:
\[
\norm{\mathcal{F}^{-1}\left(M\q\right)}_{\homos}\leq C
\norm{q}_{\homos}\qquad \text{and} \qquad
\norm{\mathcal{F}^{-1}\left(M\nabla\q\right)}_{\homos}\leq C
\norm{q}_{\homos}\, .
\]}
\end{lema}

The following lemma is fundamental to control the term $Q_I '(q)$
by the formula \eqref{estimatetwo}.
\begin{lema}
{\it \label{lemmasum} Assume $\epsi>0\,$ and $k\in\{1,2,...\}\,$.
Let us denote
\begin{align}
\label{lemma4.4a}&F_k (\tau):=\int_{\Omega_{k}
(\tau)}\e^{2\a-3+2\epsi}\int\int_{I_{k}(\eta)}
|\q(\eta-\tau'-\xi')|^2 d\s(\tau')d\s(\xi') d\s(\eta)\, ,\\
\label{lemma4.4b}&H(\tau):=\int_{\La (\tau)}\e^{2\a-3+2\epsi
}\int_{\Gae} M\q(\eta-2\tau')^2
 d\s(\tau')\, d\s(\eta)\, ,\\
\label{lemma4.4c}&G(\tau):=\int_{\widetilde{\Lambda}(\tau)}\e^{2\a-3+2\epsi}\int_{\Phi(\eta)}|M\q(\eta-\tau')|^2
d\s(\tau') d\s(\eta)\, ,
\end{align}
where
\begin{align}
\label{omegak}&\Omega_{k}(\tau):=\left\{\eta\in\Lambda (\tau) :
|\eta-\tau|\geq
C_{1}\, 2^{-k}\e\right\}\, ,\\
\label{fieta}&\Phi(\eta):=\{\ \tau'\in\,\Gae : |\tau' |\leq
1\}\quad\text{and}\quad\widetilde{\Lambda}(\tau):=\{\eta\in\Lambda(\tau)
: \e\geq 10\}\, ,
\end{align}
and $I_k (\eta)$ is given in \eqref{7b}. Then

$(i)$ If $\,\tau\in\Rdos\setminus\{0\}\,$, $ F_k (\tau)\leq C\,
\|q\|_{\dot{W}^{\a-1+\epsi,2}}^2 \, $, where $C$ only depends on
$C_1 \,$.

$(ii)$ For any $\tau\in\Rdos\,$ such that $|\tau|>1\,$, and for
any $\a \, ,\,\epsi$ so that $0<\a + \epsi<2\,$:
\[
H(\tau)\leq C\, \left(\|q\|_{\dot{W}^{\a-1+\epsi,2}}^2
+\|q\|_{L^{2}}^2 |\tau|^{2\a-2+2\epsi }\right)\, .
\]

$(iii)$ For any $\tau\in\Rdos\setminus\{0\}\,$ such that
$\taum\leq 1\,$, and provided that $0<\a+\epsi <2\,$:
\[
G(\tau)\leq C\norm{q}_{\homoame}^2 \, .
\]
}
\end{lema}

\begin{nota}
With respect to part $(i)$ of this lemma, we need that
$0<C_1<2\,$, but in fact we always apply this lemma with $C_1
<1\,$. Note that
$F_k (\tau)$ is uniformly bounded in $k\,$.\\
\end{nota}

{\it Proof of $(i)\,$.} For fixed $\tau$ we set an orthonormal
reference $\{e_{1},e_{2}\}\,$ of $\Rdos\,$, for which
$\tau=|\tau|e_{1}\,$. We write $\eta(s) = \taum e_{1} + s\,
e_{2}\,$, $s\in\R\,$. Let $h(s):=\es =(|\tau|^2 +s^2 )^{{1\over
2}}\,$. Since $|s|=|\eta(s)-\tau|\geq C_{1}2^{-k}h(s)\,$, we have
$|s|\geq C_{1} 2^{-k}\,\taum\,$. It is true that
$d\s(\eta(s))=ds\,$. We have
\[
F(\tau)= \int_{|s|\geq \, C_1 2^{-k}\taum
}(h(s))^{2\a-3+2\epsi}\int\int_{I_{k}
(\eta(s))}|\q(\eta(s)-\tau'-\xi')|^2 d\s(\tau')d\s(\xi')ds\, .
\]
Take the change of variables given by
\[
\xi'={\eta(s)\over 2}+{h(s)\over 2}\, v\quad \text{ and }\quad
\tau'={\eta(s)\over 2}+{h(s)\over 2}\, u \, ,
\]
with $u,v\in S^1 \,$. It holds that $d\s(\xi')=Ch(s) d\s(v)$ and
$d\s(\tau')=Ch(s) d\s(u)\,$. Since $|\eta(s)-\xi'-\tau'|\geq
{h(s)\over 100}\,$, hence $1+u\cdot v\geq {1\over 5000}\,$. It
holds $|u-v|\leq 4\cdot 2^{-k}\,$. We write
\begin{align*}
F(\tau )&\leq C\,\int_{|s|\geq \, C_1 2^{-k}\taum }\left(
h(s)\right)^{2\a-1+2\epsi}\int_{S^1 }\int_{A(v,\,
k)}\left|\q\left(-{h(s)\over
2}(u+v)\right)\right|^2 d\s(u) d\s(v) ds\\
&\leq C\, \suma_{j=1}^{2^k } \int_{|s|\geq \, C_1 2^{-k}\taum
}\left( h(s)\right)^{2\a-1+2\epsi}\int_{A_{j}}\int_{
\widetilde{A}_{j}}\left|\q\left(-{h(s)\over
2}(u+v)\right)\right|^2 d\s(u) d\s(v) ds\, ,
\end{align*}
where $A(v,k):=\{u\in S^1 \,:\, |u-v|\leq 4\cdot 2^{-k} \,\text{
 and  }\,\, 1+u\cdot v\geq {1\over 5000} \}\,$ and $\{A_{j} : 1\leq j\leq
2^k \,\}$ is a cover with finite overlapping of the circumference
$S^1 \,$. Each $A_{j}$ has arc-length $2^{-k}\,$. Also,
$\widetilde{A}_{j} :=\{ u\in S^1 : |u-v|\leq 4\cdot 2^{-k}\, ,\,\,
\,\, 1+u\cdot v\geq {1\over 5000}\,\text{ , for some }v\in A_{j}
\}\,$. Take the change of variables $ u(\t)=\cos\t\,
e_{1}+\sin\t\, e_{2}\, ,\; \t\in[0,2\pi)\, , \, u\in S^1 \, , $
with $d\s(u(\t))=d\t\,$. By Fubini's theorem,
\begin{align*}
F(\tau)&\leq C\, \suma_{j=1}^{2^k }\int_{A_{j}}\int_{|s|\geq \,
C_1 2^{-k}\taum }\int_{\{\t\, : \, u(\t)\,\in
\,\widetilde{A}_{j}\}} \left(
h(s)\right)^{2\a-1+2\epsi}\left|\q\left(-{h(s)\over
2}(u+v)\right)\right|^2 d\t ds d\s(v)\, .
\end{align*}
For any frozen $j\in\Z$ with $1\leq j\leq 2^k \,$ and $v\in
A_{j}\,$, take the change of variables $(s,\t )\,\ra\,
\lambda=(\l_{1}, \l_{2})\in\Rdos\,$ given by
\[
\l = -{h(s)\over 2}(u(\t)+v)= -{h(s)\over 2}((\cos\t
+v_{1})e_{1}+(\sin\t + v_{2})e_{2})\, .
\]
It holds $dsd\t = {4\over |s(1+u(\t)\cdot\, v)|}\, d\l\,$. For any
$j\geq 1\,$, we consider the proper cone
\[
H_{j}:=\left\{ r \left({u+v\over 2}\right) : r<0,\, v\in A_{j},\,
u\in\widetilde{A}_{j} \right\}\, .
\]
Since $1+u(\t)\cdot v\geq C\,$, $h(s)\sim |\l|\,$. We know that
for $|s|\geq C_1 2^{-k} \taum\,$ we also have that $|s|\geq
C(C_{1}) 2^{-k}\,h(s)\,$, $\s(A_{j})\sim 2^{-k}\,$ and the family
$\{H_{j} : 1\leq j\leq 2^k \, \}$ has finite overlap with constant
independent of $k\,$. Then
\begin{align*}
F(\tau)&\leq C\, 2^k \suma_{j=1}^{2^k
}\int_{A_{j}}\int_{H_{j}}|\l|^{2\a-2+2\epsi}|\q(\l)|^2
d\l\, d\s(v) = C \suma_{j=1}^{2^k }\int_{H_j } |\lambda |^{2\a-2+2\epsi }|\q(\lambda)|^2 d\lambda\\
&\leq C\,\int_{\Rdos}|\l|^{2\a-2+2\epsi}|\q(\l)|^2 d\l = C\,
\|q\|_{\dot{W}^{\a-1+\epsi , 2}}^2 \, .
\end{align*}

\begin{flushright} $\square $ \end{flushright}

{\it Proof of $(ii)\,$.} We follow the same lines of the previous
point but now we do not need the finite overlapping cover for $S^1
$. The variable $s$ takes real values in all the line. In the same
way, take the change $\tau' = {\eta(s)\over 2} + {|\eta(s)|\over
2}u\,$, with $u\in S^1 $ and we parametrize $u$ by $\t\in [0,2\pi
)\,$. We take the change of variables $(s,\t)\,\ra\, \l=(\l_{1},
\l_{2})$ given by
\[
\l=\eta(s)-2\tau'(\t)=-|\eta(s)| u(\t)=\, -|\eta(s)| (\cos\t\,
e_{1}+\sin\t\, e_{2})\, .
\]
Now $dsd\t = {d\l\over |s|} = {d\l\over (\, |\eta(s)|^2 - |\tau|^2
\,)^{{1\over 2}}}\,$, and $|\l|=|\eta(s)| \geq |\tau|\,$. So, we
obtain that
\begin{align*}
H(\tau)&\leq \left(\,\int_{\{\l : |\l|\geq (1+|\tau|^2 )^{{1\over
2}}\}} + \int_{\{\l : |\tau|\leq |\l|\leq (1+|\tau|^2 )^{{1\over
2}}\}}\,\right) |\l|^{2\a-2+2\epsi }M\q(\l)^2 {d\l\over (|\l|^2
-|\tau|^2
)^{{1\over 2}}}\\
&=: J_{1} + J_{2}\, .
\end{align*}
The first integral has no difficulties, indeed, $J_{1}\leq
C\norm{q}_{\homoame}^2 \,$, by lemma \ref{radial} provided that
$0<\a+\epsi<2\,$. By lemma \ref{remark}, $M\q(\l)\leq
C\norm{q}_{L^2 }$ and taking polar coordinates we get $J_{2}\leq
C\norm{q}_{L^2 }^2 \taum^{2\a -2+2\epsi}\,$.

\begin{flushright} $\square $ \end{flushright}

{\it Proof of $(iii)\,$.} Let $\tau$ be in $\Rdos\setminus
\{0\}\,$ such that $\taum\leq 1\,$. Following the same scheme as
the last point, we parametrize the variable $\eta $ by $s\in\R\,$,
take the change $\tau' = {\eta(s)\over 2}+{|\eta(s)|\over 2}\,
u\,$, with $u\in S^1 \,$, and parametrize $u$ by $\t\in
(0,2\pi]\,$. Finally take the change $(s,\t)\,\ra\,
\l=(\l_{1},\l_{2})\,$, given by $\l = \eta(s)-\tau'(\t) = {1\over
2}\left[(\taum -|\eta(s)|\cos\t ) e_{1} + (s-|\eta(s)| \sin\t )
e_{2} \right]\,$. In this case, the Jacobian for this change is
$dsd\t = {2\over |\l_{2}|}\, d\l_{1} d\l_{2}\,$. The condition
$|\tau|, |\tau'(\t)|\leq 1\,$ guarantees that the angle between
$\tau-\eta(s)\,$ and $\tau'(\t)-\eta(s)$ is uniformly bounded by
an acute angle. Remember that $|\eta(s)|\geq 10\,$. So,
$|\l_{2}|\sim |\l|\,$. That condition also implies that a positive
constant $C<1\,$ exists such that $|\l|=|\eta(s)-\tau'(\t)|\geq
C\e\,$, hence $|\eta(s)|\sim |\l|\,$. It holds
\[
G(\tau)\leq C \int_{\{\l\in\Rdos : \, |\l|>\, C'\}}
|\l|^{2\a-2+2\epsi}M\q(\l)^2 \, {d\l\over |\l|} \leq
C\,\norm{\mathcal{F}^{-1}\left(M\q\right)}_{\homoame}^2 \, ,
\]
and lemma \ref{radial} ends the proof.

\begin{flushright}
$\square $
\end{flushright}

The following lemma becomes essential to bound the terms
$Q_{j}''(q)$ in \eqref{qsegundaj}.

\begin{lema}\label{lemadelta}
{\it Let \begin{equation}
\label{expdelta}\widehat{Q''_{\d}(q)}(\eta):=\chi_{(\d^{-1},+\infty)}(\e){1\over
\e^3
}\int_{\Gad}\int_{\Gae}\left|\q(\xi)\q(\eta-\tau)\q(\tau-\xi)\right|\,
d\s(\xi)d\tau\, ,
\end{equation}
where $\Gad$ is the annulus given by
\begin{equation}
\label{gad}\Gad :=\bigg\{\tau\in\Rdos : \left|\, |\tau-{\eta\over
2}|-{\e\over 2}\,\right|\leq \d\e \bigg\}\, .
\end{equation}
Then there exist $\d_{0}\, ,\, C(\d_{0})\, ,\, \beta$ so that
$\d_{0}>0\,$, $C(\d_{0})>0\,$, $\beta>1$ and for any $\d$ with
$0<\d\leq \d_{0}\,$:
\[
\biggnorm{Q''_{\d}(q)}_{\homoa}\leq
C(\d_{0})\d^{\b}\norm{q}_{\homoame}\left(\norm{q}_{L^2
}\norm{q}_{\homomenos} + \norm{q}_{L^2 }^2 \right) ,
\]
where $\a\in\R$ and $\epsi >0$ satisfy that $0<\a+\epsi<2\,$.}
\end{lema}

We omit the proof of this lemma. We know that $d\s_{\eta
}(\tau)=\displaystyle\lim_{\d\ra 0} {1\over \d\e }\, \chi_{\Gad
}(\tau) d\tau $, where $d\s_{\eta }(\tau)$ denotes the measure on
$\Gae$ induced by $d\tau $. According to this, $\widehat{Q_{\d
}''(q)}(\eta)\sim\d\widehat{Q'(q)}(\eta)\,$, if $0<\d<\d_0 \ll
1\,$. Heuristically the estimate for $Q_{\d}''(q)$ is the one for
$Q'(q)$ multiplied by $\d\,$. We have to pay with a fraction of
derivatives in $\norm{Q'(q)}_{\homoa}$ in order to gain the factor
$\d^{\b }$ with $\b>1\,$. So, the reader must not be surprised by
the lemma whose proof follows the lines of the estimate of
$Q'(q)\,$.

\begin{notanum}
\label{notalemadelta}If we substitute \eqref{expdelta} for
\[
\widehat{Q_{\d}''}(f,g,h)(\eta)
:=\chi_{(\d^{-1},+\infty)}(\e){1\over \e^3
}\int_{\Gad}\int_{\Gae}\left|\hat{f}(\xi)\hat{g}(\eta-\tau)\hat{h}(\tau-\xi)\right|\,
d\s(\xi)d\tau\, ,
\]
with $f,g,h\in W^{s ,\, 2}\,$ and $\, -1<s<1\,$, just imitating
the proof of the control of the spherical term $Q'(q)\,$, we get
that there exist $\beta>1$ and $C(\d_{0})>0$ so that
\begin{align*}
\norm{Q_{\d}''(f,g,h)}_{\homoa}& \leq
C(\d_{0})\d^{\beta}\norm{g}_{L^2 }\big[\,\norm{f}_{L^2
}\norm{h}_{\homoame} +
\norm{f}_{\homomenos}\norm{h}_{\homoame}\\
&\qquad\qquad\qquad\qquad\qquad\qquad\qquad\qquad + \norm{h}_{L^2
}\norm{f}_{\homoame}\big] ,
\end{align*}
for any $\d$ such that $0<\d\leq\d_{0}\,$.\\
\end{notanum}

\begin{lema}
\label{lemadeltados}{\it Let $\a\in\R$ and $\epsi >0$ such that
$0<\a +\epsi <2\,$. Let $f,g,h\in W^{s ,\, 2}$ for all $s\in\R$
with $-1<s<1\,$. Let
\[
\widehat{Q_{\ast}''(q)}(\eta):=\chi(\eta)\,{1\over \e^2 }
\int_{\Gainf}\int_{\Gae}\left|\widehat{f}(\xi)\widehat{g}(\eta-\tau)\widehat{h}(\tau-\xi)\right|\,
d\s(\xi) d\tau\, ,
\]
where $\Gainf$ is the annulus given by \eqref{gainf}. Then there
exists a constant $C>0$ such that
\begin{align*}
\norm{Q_{\ast}''(q)}_{\homoa}  &\leq C\big[\,\norm{f}_{L^2
}\norm{g}_{L^2 }\norm{h}_{\homoame} +
\norm{f}_{\homomenos}\norm{g}_{L^2 }\norm{h}_{\homoame}\\
&\qquad\qquad\qquad\qquad\qquad\qquad\qquad +
\norm{f}_{\homoame}\norm{g}_{L^2 }\norm{h}_{L^2 }\big] .
\end{align*}}

\end{lema}

\begin{notanum}
Compare this lemma with remark to lemma \ref{lemadelta} when,
morally, $\d\sim\e^{-1}\,$. Consider that we work with a similar
term with an annulus with $\d\sim\e^{-1}\,$, but in the estimate
we claim the same gain of derivatives as in remark \ref{notalemadelta}.\\

\end{notanum}


\end{document}